\documentclass[11pt]{article}
\usepackage{mathrsfs}
\usepackage{amssymb,amsfonts,amsmath}
\usepackage{cite}

\def\D{\Delta}

\def\g{\gamma}
\def\l{\lambda}
\def\lr{\longrightarrow}
\def\o{\otimes}

\def\s{\sigma}

\def\v{\varepsilon}

\def\1{^{-1}}
\def\2{^{-2}}
\def\3{^{-3}}
\date{}
\newcounter{zlist}
  
\begin{document}
 \renewcommand{\baselinestretch}{1.2}
 \renewcommand{\arraystretch}{1.0}

 \title{\bf Rota-Baxter coalgebras and Rota-Baxter bialgebras}
\author{{\bf Tianshui Ma\footnote{Corresponding author}~~and Linlin Liu}\\
 {\small  $^1$Department of Mathematics, School of Mathematics and Information Science, }\\
 {\small  Henan Normal University, Xinxiang 453007, China}\\
 {\small  $^2$Henan Engineering Laboratory for Big Data Statistical Analysis and Optimal Control, }\\
 {\small School of Mathematics and Information Sciences, Henan Normal University, }\\
 {\small Xinxiang 453007, China}\\
  {\small  E-mail: matianshui@yahoo.com (T.~Ma); liulinlin2016@163.com (L. Liu)}}

 \maketitle

 \begin{center}
 \begin{minipage}{12.cm}

 {\bf \begin{center}ABSTRACT \end{center}}
 We give a construction of Rota-Baxter coalgebras from Hopf module coalgebras and also derive the structures of the pre-Lie coalgebras via Rota-Baxter coalgebras of different weight. Finally, the notion of Rota-Baxter bialgebra is introduced and some examples are provided.

 \vskip 0.1cm



 \end{minipage}
 \end{center}
 \normalsize
 \vskip0.2cm

 \section{Introduction }  \hskip\parindent
 A Rota-Baxter algebra, also called a Baxter algebra, is an associative algebra with a linear operator which generalizes the algebra of continuous functions with the integral operator. More precisely, for a given field $K$ and $\l\in K$, a Rota-Baxter $K$-algebra (of weight $\l$) is a $K$-algebra $R$ together with a $K$-linear map $P : R \lr R$ such that
 $$
 P(x)P(y)=P(xP(y))+P(P(x)y)+\l P(xy) \eqno (1)
 $$
 for all $x, y\in R$. Such a linear operator is called a Rota-Baxter operator (of weight $\l$). Rota-Baxter algebras were introduced in \cite{Ro} in the context of differential operators on commutative Banach algebras and since \cite{Baxter} intensively studied in probability and combinatorics, and more recently in mathematical physics, such as dendriform algebras (see \cite{BGN,Br,CM,GK}), free Rota-Baxter algebras (see \cite{BCQ,EFG1,EFG2,Guo2,GK1}), Lie algebras (see \cite{BGLW,LHB}), multiple zeta
 values (see \cite{EFG1,GZ}) and Connes-Kreimer renormalization theory in quantum field theory (see \cite{EFG3}), etc. One can refer to the book \cite{Guo1} for the detailed theory of Rota-Baxter algebras.

 In 2014,  Jian used Yetter-Drinfeld module algebras to construct Hopf module algebras, and then the corresponding Rota-Baxter algebras were derived (see \cite{Jian1}). Based on the dual method in the Hopf algebra theory, Jian and Zhang in \cite{JZ} defined the notion of Rota-Baxter coalgebras and also provided various examples of the new object, including constructions by group-like elements and by smash coproduct.

 In this paper, we will continue to investigate the properties of Rota-Baxter coalgebras. We give a construction of Rota-Baxter coalgebras from Hopf module coalgebras and also derive the structures of the pre-Lie coalgebras via Rota-Baxter coalgebras of different weight. Finally, the notion of Rota-Baxter bialgebra is introduced  and some examples are provided.

 \section{Preliminaries}  \hskip\parindent
 Throughout this paper, we assume that all vector spaces, algebras, coalgebras and tensor products are defined over a field $K$. An algebra is always assumed to be associative, but not necessarily unital. A coalgebra is always assumed to be coassociative, but not necessarily counital. Now, let $C$ be a coalgebra. We use Sweedler's notation for the comultiplication (see \cite{Mon}): $\Delta(c)=c_{1}\o c_{2}$ for any $c\in C$. Denote the category of left $H$-comodules by $^H \hbox{Mod}$. For $(M, \rho_L)\in$ $^H \hbox{Mod}$, write: $\rho_L(x)=x_{(-1)}\o x_{(0)} \in H\o M$, for all $x \in M$. Denote the category of right $H$-comodules by $ \hbox{Mod}^H$. For $(M, \rho_R)\in$ $ \hbox{Mod}^H$, write: $ \rho_R(x)=x_{(0)}\o x_{(1)} \in M\o H$, for all $x \in M$.  Given a $K$-space $M$, we write $id_M$ for the identity map on $M$.

 In what follows, we recall  some useful definitions and results which will be used later (see \cite{Guo1,Jian1,JZ,Mon}).

 {\bf Hopf module} (see \cite{Mon}): A right (resp. left) $H$-Hopf module is a vector space $M$ equipped simultaneously with a right (resp. left) $H$-module structure and a right (resp. left) $H$-comodule structure $\rho_{R}$ (resp. $\rho_{L}$) such that
 $$
 \rho_{R}(m\cdot h)=\rho_{R}(m)\D(h).    \eqno(2)
 $$
 (resp. $\rho_{L}(h\cdot m)=\D(h)\rho_{L}(m)$) for all $h\in H$ and $m\in M$.

 {\bf Remark} Let $M$ be a right $H$-Hopf module. We denote by $M^{R}$ the subspace of right coinvariants, i.e.,
 $$
 M^{R}=\{m\in M\mid \rho_{R}(m)=m\o 1_H\}.
 $$
 Then the map $P_{R} : M \lr M$ given by $P_{R}(m)=m_{(0)}\cdot S(m_{(1)})$ for $m\in M$ is the projection from $M$ onto $M^{R}$. Similarly, for a left $H$-Hopf module, we can define the subspace $M^{L}$ of left coinvariants, i.e.,
 $$
 M^{L}=\{m\in M\mid \rho_{L}(m)=1_H\o m\},
 $$
 and have the projection formula
 $$
 P_{L}(m)=S(m_{(-1)})\cdot m_{(0)}.
 $$

 \smallskip
 {\bf Hopf module algebra} (see \cite{Jian1}):  A right $H$-Hopf module algebra is a right $H$-Hopf module $M$ together with an associative multiplication $\mu : M \o M\lr M$ (write $\mu(m\o m')=mm'$) such that for any $h\in H$ and $m, m'\in M$,
 $$
 (mm')\cdot h= m(m'\cdot h),        \eqno(3)
 $$
 $$
 \rho_{R}(mm')=m_{(0)}m'_{(0)}\o m_{(1)}m'_{(1)}.  \eqno(4)
 $$

 {\bf Remarks} (1) The compatibility conditions between $\mu$ and the Hopf module structure mean that  $\mu$ is both a module morphism and a comodule morphism,where the right $H$-module action and right $H$-comodule coaction on $M\o M$ are given by
 $$
 M\o M\o H\lr M\o M ~~~(m\o m'\o h\mapsto m\o m'\cdot h),
 $$
 $$
 M\o M\lr M\o M\o H ~~~(m\o m'\mapsto m_{(0)}\o m'_{(0)}\o m_{(1)}m'_{(1)})
 $$
 for all $h\in H$ and $m, m'\in M$.\newline
 \indent{\phantom{\bf Remarks}} (2) We can define a left $H$-Hopf module algebra by satisfying the compatibility conditions $h\cdot (mm')=(h\cdot m)m'$ and $\rho_{L}(mm')=m_{(-1)}m'_{(-1)}\o m_{(0)}m'_{(0)}$.

 \smallskip

 {\bf Rota-Baxter coalgebra} (see \cite{JZ}): Let $\g$ be an element in $K$. A pair $(C, Q)$ is called a Rota-Baxter coalgebra of weight $\gamma$ if $C$ is a coalgebra and $Q$ is a linear endomorphism of $C$ satisfying that
 $$
 (Q\o Q)\D=(id\o Q)\D Q+(Q\o id)\D Q+\gamma\D Q,
 $$
 that is, for all $c\in C$,
 $$
 Q(c_1)\o Q(c_2)=Q(c)_1\o Q(Q(c)_2)+Q(Q(c)_1)\o Q(c)_2+ \gamma Q(c)_1\o Q(c)_2. \eqno(5)
 $$
 The map $Q$ is called a Rota-Baxter operator.

 A Rota-Baxter coalgebra is called idempotent if the Rota-Baxter operator $Q$ satisfies $Q^{2}=Q$.

 \smallskip
 {\bf Yetter-Drinfeld category} (see \cite{Mon}): A left~$H$-Yetter-Drinfeld~module is a vector space $V$ equipped with a left $H$-module structure $\cdot$ and a left $H$-comodule structure $\rho_L$ such that for all $h\in H$ and $v\in V$,
 $$
 h_1v_{(-1)}\o h_2\cdot v_{(0)}=(h_1\cdot v)_{(-1)}h_2\o (h_1\cdot v)_{(0)}.\eqno(6)
 $$
 The category of left $H$-Yetter-Drinfeld modules, denoted by $^{H}_{H}\mathcal{YD}$, consists of the following data: Its objects are left $H$-Yetter-Drinfeld modules and morphisms are linear maps which are both module and comodule morphisms.

 {\bf Remarks} (1) A coalgebra $(C,\D)$ in $^{H}_{H}\mathcal{YD}$ is a coassociative coalgebra such that the underlying space $C$ is an object in $^{H}_{H}\mathcal{YD}$ and $(C,\D)$ is both a module-coalgebra and a comodule-coalgebra. More precisely, for all $h\in H$ and $c\in C$ we have
 $$
 \D(h\cdot c)=h_1\cdot c_1\o h_2\cdot c_2,\eqno(7)
 $$
 $$
 c_{(-1)}\o c_{(0)1}\o c_{(0)2}=c_{1(-1)}c_{2(-1)}\o c_{1(0)}\o c_{2(0)}.\eqno(8)
 $$ \newline
 \indent{\phantom{\bf Remarks}} (2)
 Let $(C,\D)$ be a coalgebra in $^{H}_{H}\mathcal{YD}$. The smash coproduct coalgebra $C\times H$ is defined to be $C\o H$ as a vector space and
 $$
 \D(c\o h)=c_1\o c_{2(-1)}h_1\o c_{2(0)}\o h_2.\eqno(9)
 $$

 \smallskip

 \section{Construction of Rota-Baxter coalgebras}  \hskip\parindent
 In this section, we give a construction of Rota-Baxter coalgebras via Hopf module coalgebras.

 {\bf Definition 3.1} A right $H$-Hopf module coalgebra is a right $H$-Hopf module $M$ together with coassociative comultiplication $\D : M \lr M\o M$ ($\D(m)=m_1\o m_2$) such that for any $h\in H$ and $m\in M$,
 $$
 m_{(0)1}\o m_{(0)2}\o m_{(1)}=m_1\o m_{2(0)}\o m_{2(1)}, \eqno(10)
 $$
 $$
 \D(m\cdot h)=m_1\cdot h_1\o m_2\cdot h_2.   \eqno(11)
 $$

 {\bf Remarks} (1) The compatibility conditions between $\D$ and the Hopf module structure mean that  $\D$ is both a module morphism and a comodule morphism, where the right $H$-module action and right $H$-comodule coaction on $M\o M$ are given by
 $$
 M\o M \o H\lr M\o M ~~~(m\o m'\o h\mapsto m\cdot h_1\o m'\cdot h_2),
 $$
 $$
 M\o M\lr M\o M\o H ~~~(m\o m'\mapsto m\o m'_{(0)}\o m'_{(1)})
 $$
 for all $h\in H$ and $m, m'\in M$. \newline
 \indent{\phantom{\bf Remarks}} (2) We can define a left $H$-Hopf module coalgebra by satisfying the compatibility conditions:
 $$
 m_{(-1)}\o m_{(0)1}\o m_{(0)2}=m_{1(-1)}\o m_{1(0)}\o m_2,
 $$
 and
 $$
 \D(h\cdot m)=h_1\cdot m_1\o h_2\cdot m_2,
 $$
 where $h\in H$ and $m\in M$.

 \smallskip

 {\bf Theorem 3.2} Let $(C,\D)$ be a right $H$-Hopf module coalgebra. Then $(C, P_{R})$ is a Rota-Baxter coalgebra of weight $-1$.

 {\bf Proof.} Since $C$ is a right $H$-comodule, for all $c\in C$ we have
 $$
 c_{(0)}\o c_{(1)1}\o c_{(1)2}=c_{(0)(0)}\o c_{(0)(1)}\o c_{(1)}.\eqno(12)
 $$
 We check the condition $(5)$ for $P_{R}$ as follows:
 \begin{eqnarray*}
 &&((id\o P_{R})\D P_{R}+(P_{R}\o id)\D P_{R}-\D P_{R})(c)\\
 &=&(id\o P_{R})\D P_{R}(c)+(P_{R}\o id)\D P_{R}(c)-\D P_{R}(c)\\
 &=&(id\o P_{R})\D(c_{(0)}\cdot S(c_{(1)}))+(P_{R}\o id)\D(c_{(0)}\cdot S(c_{(1)}))-\D(c_{(0)}\cdot S(c_{(1)}))\\
 &=&(c_{(0)}\cdot S(c_{(1)}))_1\o P_{R}((c_{(0)}\cdot S(c_{(1)}))_2)+P_{R}((c_{(0)}\cdot S(c_{(1)}))_1)\o (c_{(0)}\cdot S(c_{(1)}))_2\\
 &&-(c_{(0)}\cdot S(c_{(1)}))_1\o(c_{(0)}\cdot S(c_{(1)}))_2\\
 &\stackrel{(11)}{=}&c_{(0)1}\cdot S(c_{(1)2})\o (c_{(0)2}\cdot S(c_{(1)1}))_{(0)}\cdot S((c_{(0)2}\cdot S(c_{(1)1}))_{(1)})+(c_{(0)1}\cdot S(c_{(1)2}))_{(0)}\\
 &&\cdot S((c_{(0)1}\cdot S(c_{(1)2}))_{(1)})\o c_{(0)2}\cdot S(c_{(1)1})-c_{(0)1}\cdot S(c_{(1)2})\o c_{(0)2}\cdot S(c_{(1)1})\\
 &\stackrel{(2)}{=}&c_{(0)1}\cdot S(c_{(1)2})\o(c_{(0)2(0)}\cdot S(c_{(1)1})_{1})\cdot S(c_{(0)2(1)}S(c_{(1)1})_{2})+c_{(0)1(0)}\cdot S(c_{(1)2})_{1}\\
 &&\cdot S(c_{(0)1(1)}S(c_{(1)2})_{2})\o c_{(0)2}\cdot S(c_{(1)1})-c_{(0)1}\cdot S(c_{(1)2})\o c_{(0)2}\cdot S(c_{(1)1})\\
 &=&c_{(0)1}\cdot S(c_{(1)2})\o c_{(0)2(0)}\cdot ((S(c_{(1)1})_{1}S(S(c_{(1)1})_{2}))S(c_{(0)2(1)}))+c_{(0)1(0)}\cdot ((S(c_{(1)2})_{1}\\
 &&S(S(c_{(1)2})_{2}))S(c_{(0)1(1)}))\o c_{(0)2}\cdot S(c_{(1)1})-c_{(0)1}\cdot S(c_{(1)2})\o c_{(0)2}\cdot S(c_{(1)1})\\
 &=&c_{(0)1}\cdot S(c_{(1)2})\o c_{(0)2(0)}\cdot \v(S(c_{(1)1}))S(c_{(0)2(1)})+c_{(0)1(0)}\cdot S(c_{(0)1(1)})\\
 &&\o c_{(0)2}\cdot S(c_{(1)1})\v(S(c_{(1)2}))-c_{(0)1}\cdot S(c_{(1)2})\o c_{(0)2}\cdot S(c_{(1)1})\\
 &=&c_{(0)1}\cdot S(c_{(1)})\o c_{(0)2(0)}\cdot S(c_{(0)2(1)})+c_{(0)1(0)}\cdot S(c_{(0)1(1)})\o c_{(0)2}\cdot S(c_{(1)})\\
 &&-c_{(0)1}\cdot S(c_{(1)2})\o c_{(0)2}\cdot S(c_{(1)1})\\
 &\stackrel{(10)}{=}&c_{(0)(0)1}\cdot S(c_{(1)})\o c_{(0)(0)2}\cdot S(c_{(0)(1)})+c_{(0)1(0)}\cdot S(c_{(0)1(1)})\o c_{(0)2}\cdot S(c_{(1)})\\
 &&-c_{(0)1}\cdot S(c_{(1)2})\o c_{(0)2}\cdot S(c_{(1)1})
  \end{eqnarray*}
   \begin{eqnarray*}
 &\stackrel{(12)}{=}&c_{(0)1}\cdot S(c_{(1)2})\o c_{(0)2}\cdot S(c_{(1)1})+c_{(0)1(0)}\cdot S(c_{(0)1(1)})\o c_{(0)2}\cdot S(c_{(1)})\\
 &&-c_{(0)1}\cdot S(c_{(1)2})\o c_{(0)2}\cdot S(c_{(1)1})\\
 &=&c_{(0)1(0)}\cdot S(c_{(0)1(1)})\o c_{(0)2}\cdot S(c_{(1)})\\
 &\stackrel{(10)}{=}&c_{1(0)}\cdot S(c_{1(1)})\o c_{2(0)}\cdot S(c_{2(1)})=(P_{R}\o P_{R})\D(c),
 \end{eqnarray*}
 finishing the proof.                                                                      \hfill $\square$

  {\bf Remark} Let $(C,\D)$ be a left $H$-Hopf module coalgebra. Then $(C, P_{L})$ is a Rota-Baxter coalgebra of weight $-1$.

 \smallskip

 {\bf Example 3.3} Let $C$ be a bialgebra and $H$ a Hopf algebra with the antipode $S$. Suppose that there are two bialgebra maps: $i: H\lr C$ and $\pi: C\lr H$ such that $\pi\circ i=id_{H}$, i.e., $C$ is a bialgebra with a projection (see \cite{Ra}). Set $\Pi=id_{C}\star (i\circ S\circ \pi)$, where $\star$ is the convolution product on End($C$).
The right $H$-module structure and right $H$-comodule structure on $C$ are given by the following:
 $$
 c\cdot h=ci(h),
 $$
 $$
 \rho_{R}(c)=c_1\o \pi(c_2),
 $$
 for all $c\in C$ and $h\in H$. Then $(C, \Pi)$ is a Rota-Baxter coalgebra of weight $-1$.

 {\bf Proof.} For all $h\in H$ and $c\in C$, firstly $C$ is a right $H$-Hopf module  with the above module structure and comodule structure, since
 \begin{eqnarray*}
 \rho_{R}(c\cdot h)&=&\rho_{R}(ci(h))=(ci(h))_1\o \pi((ci(h))_2)\\
 &=&c_1i(h_1)\o \pi(c_2)h_2=c_1\cdot h_1\o \pi(c_2)h_2=\rho_{R}(c)\D(h).
 \end{eqnarray*}

 Secondly, we have
  \begin{eqnarray*}
  (\D\o id)\rho_{R}(c)
  &=&(\D\o id)(c_1\o \pi(c_2))=c_{11}\o c_{12}\o \pi(c_{2})\\
  &=&c_1\o c_{21}\o \pi(c_{22})=(id\o \rho_{R})\D(c)
  \end{eqnarray*}
  and $\D(c\cdot h)=c_1i(h_1)\o c_2i(h_2)=c_1\cdot h_1\o c_2\cdot h_2$.  Thus $(C, \D)$ is a right $H$-Hopf module coalgebra.

 Finally, by \cite[Theorem 3]{Ra}, $\Pi$ is the projection from C to $C^{R}$.  Therefore, we get the conclusion by Theorem 3.2.                                    \hfill $\square$

 \smallskip

 Next we investigate the Hopf module coalgebra structure on the smash coproduct coalgebra.

 {\bf Theorem 3.4} Let $(C,\D)$ be a coalgebra in $^{H}_{H}\mathcal{YD}$. Then: (1) The smash coproduct coalgebra $C\times H$ is a right $H$-Hopf module coalgebra with the following right $H$-Hopf module structure: for all $h, x\in H$ and $c\in C$,
 $$
 (c\o h)\cdot x=c\o hx,
 $$
 $$
 \rho_{R}(c\o h)=(c\o h_1)\o h_2.
 $$
 Hence the map $P_{R}$ given by
 $$
 P_{R}(c\o h)=c\o \v(h)1_H
 $$
 is an idempotent Rota-Baxter operator of weight $-1$.

 (2) The smash coproduct coalgebra $C\times H$ is a left $H$-Hopf module coalgebra with the following left $H$-Hopf module structure:
 $$
 x\cdot(c\o h)=x_1\cdot c\o x_2h,
 $$
 $$
 \rho_{L}(c\o h)=c_{(-1)}h_1\o (c_{(0)}\o h_2)
 $$
 Hence the map $P_{L}$ given by
 $$
 P_{L}(c\o h)=S(c_{(-1)2}h_2)\cdot c_{(0)}\o S(c_{(-1)1}h_1)h_3
 $$
 is an idempotent Rota-Baxter operator of weight $-1$.

 {\bf Proof.} In what follows, we only prove the compatibility conditions, and the constructions of the right and left $H$-Hopf module structures given above are obvious.

 (1) For all $c\in C$ and $h,x\in H$, we have
 \begin{eqnarray*}
 \D((c\o h)\cdot x)
 &=&\D(c\o hx)\stackrel{(9)}{=}c_1\o c_{2(-1)}(hx)_1\o c_{2(0)}\o (hx)_2\\
 &=&c_1\o (c_{2(-1)}h_1)x_1\o c_{2(0)}\o h_2x_2\\
 &=&(c_1\o c_{2(-1)}h_1)\cdot x_1\o (c_{2(0)}\o h_2)\cdot x_2\\
 &=&(c\o h)_1\cdot x_1\o (c\o h)_2\cdot x_2.
 \end{eqnarray*}
 and
 \begin{eqnarray*}
 (id\o \rho_{R})\D(c\o h)
 &\stackrel{(9)}{=}&(id\o \rho_{R})(c_1\o c_{2(-1)}h_1\o c_{2(0)}\o h_2)\\
 &=&c_1\o c_{2(-1)}h_1\o c_{2(0)}\o h_{21}\o h_{22}\\
 &=&c_1\o c_{2(-1)}h_{11}\o c_{2(0)}\o h_{12}\o h_{2}\\
 &=&(\D\o id)(c\o h_1\o h_2)=(\D\o id)\rho_{R}(c\o h).
 \end{eqnarray*}

 (2) Since $C$ is a left $H$-comodule, for all $c\in C$ we have
 $$
 c_{(-1)}\o c_{(0)(-1)}\o c_{(0)(0)}=c_{(-1)1}\o c_{(-1)2}\o c_{(0)}.                    \eqno(13)
 $$
 While
 \begin{eqnarray*}
 &&\D(x\cdot(c\o h))=\D(x_1\cdot c\o x_2h)\\
 &&~~~~~\stackrel{(9)}{=}(x_1\cdot c)_1\o (x_1\cdot c)_{2(-1)}(x_2h)_1\o (x_1\cdot c)_{2(0)}\o(x_2h)_2\\
 &&~~~~~\stackrel{(7)}{=}x_{11}\cdot c_1\o(x_{12}\cdot c_2)_{(-1)}(x_{21}h_1)\o(x_{12}\cdot c_2)_{(0)}\o x_{22}h_2\\
 &&~~~~~~=x_{1}\cdot c_1\o((x_{2}\cdot c_2)_{(-1)}x_{3})h_1\o (x_{2}\cdot c_2)_{(0)}\o x_{4}h_2\\
 &&~~~~~\stackrel{(6)}{=}x_{1}\cdot c_1\o(x_{2}c_{2(-1)})h_1\o x_{3}\cdot c_{2(0)}\o x_{4}h_2\\
 &&~~~~~~=x_{11}\cdot c_1\o x_{12}(c_{2(-1)}h_1)\o x_{21}\cdot c_{2(0)}\o x_{22}h_2\\
 &&~~~~~~=x_{1}\cdot (c_1\o c_{2(-1)}h_1)\o x_{2}\cdot (c_{2(0)}\o h_2)\\
 &&~~~~~~=x_1\cdot(c\o h)_1\o x_2\cdot(c\o h)_2,
 \end{eqnarray*}
 and
 \begin{eqnarray*}
 &&(id\o \D)\rho_{L}(c\o h)=(c\o h)_{(-1)}\o(c\o h)_{(0)1}\o(c\o h)_{(0)2}\\
 &&~~~~~~=c_{(-1)}h_1\o(c_{(0)}\o h_2)_1\o(c_{(0)}\o h_2)_2\\
 &&~~~~~~=c_{(-1)}h_1\o c_{(0)1}\o c_{(0)2(-1)}h_{21}\o c_{(0)2(0)}\o h_{22}\\
 &&~~~~~~=c_{(-1)}h_1\o c_{(0)1}\o c_{(0)2(-1)}h_{2}\o c_{(0)2(0)}\o h_{3}\\
 &&~~~~~\stackrel{(8)}{=}(c_{1(-1)}c_{2(-1)})h_1\o c_{1(0)}\o c_{2(0)(-1)}h_{2}\o c_{2(0)(0)}\o h_{3}\\
 &&~~~~~\stackrel{(13)}{=}(c_{1(-1)}c_{2(-1)1})h_1\o c_{1(0)}\o c_{2(-1)2}h_{2}\o c_{2(0)}\o h_{3}\\
 &&~~~~~~=c_{1(-1)}(c_{2(-1)}h_1)_1\o c_{1(0)}\o (c_{2(-1)}h_{1})_2\o c_{2(0)}\o h_{2}\\
 &&~~~~~~=(c_1\o c_{2(-1)}h_1)_{(-1)}\o(c_1\o c_{2(-1)}h_1)_{(0)}\o c_{2(0)}\o h_{2}\\
 &&~~~~~\stackrel{(9)}{=}(c\o h)_{1(-1)}\o (c\o h)_{1(0)}\o (c\o h)_{2}\\
 &&~~~~~~=(\rho_{L}\o id)\D(c\o h),
 \end{eqnarray*}
 thus we complete the proof.                                                                                                           \hfill $\square$

 \smallskip

 {\bf Example 3.5} The coalgebra $H$ is a coalgebra in $^{H}_{H}\mathcal{YD}$ with the following left $H$-Yetter-Drinfeld module structure:
 $$
 x\cdot h=xh,
 $$
 $$
 \rho(h)=h_1S(h_3)\o h_2.
 $$
 Then the smash coproduct on $H\o H$ is given by
 $$
 \D(h\o h')=h_1\o (h_2S(h_4))h'_1\o h_3\o h'_2,
 $$
 the corresponding Rota-Baxter operators $P_{R}$ and $P_{L}$ are given by
 $$
 P_{R}(h\o h')=h\o 1_H\v(h')
 $$
 and
 \begin{eqnarray*}
 P_{L}(h\o h')
 &=&S((h_1S(h_3))_2h'_2)\cdot h_2\o S((h_1S(h_3))_1h'_1)h'_3\\
 &=&S((h_{12}S(h_{31}))h'_2)\cdot h_2\o S((h_{11}S(h_{32}))h'_1)h'_3\\
 &=&S((h_{2}S(h_{4}))h'_2)h_2\o S((h_{1}S(h_{5}))h'_1)h'_3.
 \end{eqnarray*}

 We recall from \cite{Mon} that a coquasitriangular (or braided) Hopf algebra is a pair $(H,\s)$, where $H$ is a Hopf algebra and $\s \in (H \otimes H)^*$ such that
 \begin{eqnarray*}
 &&(BR1)\quad ~\s(1,h)=\s(h,1)=\v(h), \\
 &&(BR2)\quad ~\s(h h', h'')=\s(h, h''_{1})\s(h', h''_{2}), \\
 &&(BR3)\quad ~\s(h, h' h'')=\s(h_{1}, h'')\s(h_{2}, h'), \\
 &&(BR4)\quad ~h'_{1}h_{1}\s(h_{2}, h'_{2})=\s(h_{1},  h'_{1})h_{2}h'_{2}
 \end{eqnarray*}
 for any $h, h', h'' \in H$.

 {\bf Example 3.6} Let $(H, \sigma)$ be a coquasitriangular Hopf algebra. For any $H$-comodule $(M, \rho)$, we define $\cdot: H\o M\lr M$ by
 $ h\cdot m=\sigma(m_{(-1)}, h)m_{(0)}$.  Then $(M,\cdot, \rho)$ is a left $H$-Yetter-Drinfeld module. By \cite{Mon}, we know that any comodule-coalgebra $C$ over $H$ is moreover a coalgebra in $^{H}_{H}\mathcal{YD}$. Then Rota-baxter operators $P_{R}$ and $P_{L}$ on $C\times H$ are given by
 $$
 P_{R}(c\o h)=c\o\v(h)1_H
 $$
 and
 \begin{eqnarray*}
 P_{L}(c\o h)
 &=&S(c_{(-1)2}h_2)\cdot c_{(0)}\o S(c_{(-1)1}h_1)h_3\\
 &=&\sigma(c_{(0)(-1)}, S(c_{(-1)2}h_2))c_{(0)(0)}\o S(c_{(-1)1}h_1)h_3\\
 &\stackrel{(13)}{=}&\sigma(c_{(-1)2}, S(c_{(-1)12}h_2))c_{(0)}\o S(c_{(-1)11}h_1)h_3\\
 &=&\sigma(c_{(-1)3}, S(c_{(-1)2}h_2))c_{(0)}\o S(c_{(-1)1}h_1)h_3\\
 &=&\sigma(c_{(-1)3}, S(h_2)S(c_{(-1)2}))c_{(0)}\o S(c_{(-1)1}h_1)h_3\\
 &\stackrel{(BR3)}{=}&\sigma(c_{(-1)3}, S(c_{(-1)2}))\sigma(c_{(-1)4}, S(h_2))c_{(0)}\o S(c_{(-1)1}h_1)h_3
 \end{eqnarray*}
  respectively.

 \section{From Rota-Baxter coalgebra to pre-Lie coalgebra}  \hskip\parindent
 In this section, we construct pre-Lie coalgebras from Rota-Baxter coalgebras. First, we give the definition of pre-Lie coalgebra which is dual to the definition of pre-Lie algebra in \cite{AB}.

 {\bf Definition 4.1} A pre-Lie coalgebra is $(C,\D)$ consisting of a linear space $C$, a linear map $\D: C\lr C\o C$ and satisfying
 $$
 \D_{C}-\Phi_{(12)}\D_{C}=0,            \eqno(14)
 $$
 where $\D_{C}=(\D\o id)\D-(id\o\D)\D$ and $\Phi_{(12)}(c_1\o c_2\o c_3)=c_2\o c_1\o c_3$.

 \smallskip

 {\bf Theorem 4.2} Let $(C, \D, Q)$ be a Rota-Baxter coalgebra of weight $-1$. Define the operation $\widetilde{\D}$ on $C$ by
 $$
 \widetilde{\D}(c)=Q(c_1)\o c_2-Q(c_2)\o c_1-c_1\o c_2.
 $$
 Then $\widetilde{C}=(C,\widetilde{\D})$ is a pre-Lie coalgebra.

 {\bf Proof.} For all $c\in C$, we have
 \begin{eqnarray*}
 &&\widetilde{\D}_{C}(c)=(\widetilde{\D}\o id)\widetilde{\D}(c)-(id\o \widetilde{\D})\widetilde{\D}(c)\\
 &&~~~~~~=(\widetilde{\D}\o id)(Q(c_1)\o c_2-Q(c_2)\o c_1-c_1\o c_2)\\
 &&~~~~~~~~~~ -(id\o \widetilde{\D})(Q(c_1)\o c_2-Q(c_2)\o c_1-c_1\o c_2)\\
 &&~~~~~~=\widetilde{\D}(Q(c_1))\o c_2-\widetilde{\D}(Q(c_2))\o c_1-\widetilde{\D}(c_1)\o c_2-Q(c_1)\o\widetilde{\D}(c_2)\\
 &&~~~~~~~~~~+Q(c_2)\o\widetilde{\D}(c_1)+c_1\o\widetilde{\D}(c_2)\\
 &&~~~~~~=Q(Q(c_1)_1)\o Q(c_1)_2\o c_2-Q(Q(c_1)_2)\o Q(c_1)_1\o c_2-Q(c_1)_1\o Q(c_1)_2\o c_2\\
 &&~~~~~~~~~~ -Q(Q(c_2)_1)\o Q(c_2)_2\o c_1+Q(Q(c_2)_2)\o Q(c_2)_1\o c_1+Q(c_2)_1\o Q(c_2)_2\o c_1\\
 &&~~~~~~~~~~ -Q(c_{11})\o c_{12}\o c_{2}+ Q(c_{12})\o c_{11}\o c_{2}+ c_{11}\o c_{12}\o c_{2}-Q(c_{1})\o Q(c_{21})\o c_{22}\\
 &&~~~~~~~~~~+ Q(c_{1})\o Q(c_{22})\o c_{21}+ Q(c_{1})\o c_{21}\o c_{22}+Q(c_{2})\o Q(c_{11})\o c_{12}\\
 &&~~~~~~~~~~ -Q(c_{2})\o Q(c_{12})\o c_{11}- Q(c_{2})\o c_{11}\o c_{12}+c_{1}\o Q(c_{21})\o c_{22}\\
 &&~~~~~~~~~~ - c_1\o Q(c_{22})\o c_{21} - c_{1}\o c_{21}\o c_{22}\\
 &&~~~~~~=Q(Q(c_1)_1)\o Q(c_1)_2\o c_2-Q(Q(c_1)_2)\o Q(c_1)_1\o c_2-Q(c_1)_1\o Q(c_1)_2\o c_2\\
 &&~~~~~~~~~~ -Q(Q(c_2)_1)\o Q(c_2)_2\o c_1+Q(Q(c_2)_2)\o Q(c_2)_1\o c_1+Q(c_2)_1\o Q(c_2)_2\o c_1\\
 &&~~~~~~~~~~ -Q(c_{1})\o c_{2}\o c_{3}+ Q(c_{2})\o c_{1}\o c_{3}+ c_{1}\o c_{2}\o c_{3} -Q(c_{1})\o Q(c_{2})\o c_{3}\\
 &&~~~~~~~~~~+ Q(c_{1})\o Q(c_{3})\o c_{2}+ Q(c_{1})\o c_{2}\o c_{3}+Q(c_{3})\o Q(c_{1})\o c_{2}\\
 &&~~~~~~~~~~ - Q(c_{3})\o Q(c_{2})\o c_{1}- Q(c_{3})\o c_{1}\o c_{2}+c_{1}\o Q(c_{2})\o c_{3}\\
 &&~~~~~~~~~~ - c_1\o Q(c_{3})\o c_{2}- c_{1}\o c_{2}\o c_{3}\\
 &&~~~~~~=Q(Q(c_1)_1)\o Q(c_1)_2\o c_2-Q(Q(c_1)_2)\o Q(c_1)_1\o c_2-Q(c_1)_1\o Q(c_1)_2\o c_2\\
 &&~~~~~~~~~~ -Q(Q(c_2)_1)\o Q(c_2)_2\o c_1+Q(Q(c_2)_2)\o Q(c_2)_1\o c_1+Q(c_2)_1\o Q(c_2)_2\o c_1\\
 &&~~~~~~~~~~ +Q(c_{2})\o c_{1}\o c_{3}-Q(c_{1})\o Q(c_{2})\o c_{3}+ Q(c_{1})\o Q(c_{3})\o c_{2}\\
 &&~~~~~~~~~~ +Q(c_{3})\o Q(c_{1})\o c_{2}-Q(c_{3})\o Q(c_{2})\o c_{1}- Q(c_{3})\o c_{1}\o c_{2}+c_{1}\o Q(c_{2})\o c_{3}\\
 &&~~~~~~~~~~ - c_1\o Q(c_{3})\o c_{2}.
 \end{eqnarray*}

 Then we can verify that the condition (14) holds by using the above equality as follows.
 \begin{eqnarray*}
 &&(\widetilde{\D}_{C}-\Phi_{(12)}\widetilde{\D}_{C})(c)\\
 &&~~~~~~=Q(Q(c_1)_1)\o Q(c_1)_2\o c_2-Q(Q(c_1)_2)\o Q(c_1)_1\o c_2-Q(c_1)_1\o Q(c_1)_2\o c_2\\
 &&~~~~~~~~~~ -Q(Q(c_2)_1)\o Q(c_2)_2\o c_1+Q(Q(c_2)_2)\o Q(c_2)_1\o c_1+Q(c_2)_1\o Q(c_2)_2\o c_1\\
 &&~~~~~~~~~~ + Q(c_{2})\o c_{1}\o c_{3}-Q(c_{1})\o Q(c_{2})\o c_{3}+ Q(c_{1})\o Q(c_{3})\o c_{2}\\
 &&~~~~~~~~~~ +Q(c_{3})\o Q(c_{1})\o c_{2}-Q(c_{3})\o Q(c_{2})\o c_{1}- Q(c_{3})\o c_{1}\o c_{2}+c_{1}\o Q(c_{2})\o c_{3}\\
 &&~~~~~~~~~~ - c_1\o Q(c_{3})\o c_{2}- Q(c_1)_2\o Q(Q(c_1)_1)\o c_2+ Q(c_1)_1\o Q(Q(c_1)_2)\o c_2\\
 &&~~~~~~~~~~ + Q(c_1)_2\o Q(c_1)_1\o c_2+ Q(c_2)_2\o Q(Q(c_2)_1)\o c_1- Q(c_2)_1\o Q(Q(c_2)_2)\o c_1\\
 &&~~~~~~~~~~ - Q(c_2)_2\o Q(c_2)_1\o c_1- c_{1}\o Q(c_{2})\o c_{3}+ Q(c_{2})\o Q(c_{1})\o c_{3}\\
 &&~~~~~~~~~~ - Q(c_{3})\o Q(c_{1})\o c_{2}- Q(c_{1})\o Q(c_{3})\o c_{2}+ Q(c_{2})\o Q(c_{3})\o c_{1}\\
 &&~~~~~~~~~~ + c_{1}\o Q(c_{3})\o c_{2}- Q(c_{2})\o c_{1}\o c_{3}+ Q(c_{3})\o c_1\o c_{2}\\
 &&~~~~~~=Q(Q(c_1)_1)\o Q(c_1)_2\o c_2-Q(Q(c_1)_2)\o Q(c_1)_1\o c_2-Q(c_1)_1\o Q(c_1)_2\o c_2\\
 &&~~~~~~~~~~ -Q(Q(c_2)_1)\o Q(c_2)_2\o c_1+Q(Q(c_2)_2)\o Q(c_2)_1\o c_1+Q(c_2)_1\o Q(c_2)_2\o c_1\\
 &&~~~~~~~~~~ -Q(c_{1})\o Q(c_{2})\o c_{3}- Q(c_{3})\o Q(c_{2})\o c_{1} - Q(c_1)_2\o Q(Q(c_1)_1)\o c_2\\
 &&~~~~~~~~~~ + Q(c_1)_1\o Q(Q(c_1)_2)\o c_2+ Q(c_1)_2\o Q(c_1)_1\o c_2 + Q(c_2)_2\o Q(Q(c_2)_1)\o c_1\\
 &&~~~~~~~~~~ - Q(c_2)_1\o Q(Q(c_2)_2)\o c_1- Q(c_2)_2\o Q(c_2)_1\o c_1 + Q(c_{2})\o Q(c_{1})\o c_{3}\\
 &&~~~~~~~~~~ + Q(c_{2})\o Q(c_{3})\o c_{1}\\
 &&~~~~~~\stackrel{(5)}{=}Q(c_{11})\o Q(c_{12})\o c_2- Q(c_{12})\o Q(c_{11})\o c_2- Q(c_{21})\o Q(c_{22})\o c_1\\
 &&~~~~~~~~~~ + Q(c_{22})\o Q(c_{21})\o c_1 -Q(c_{1})\o Q(c_{2})\o c_3- Q(c_{3})\o Q(c_{2})\o c_1\\
 &&~~~~~~~~~~ + Q(c_{2})\o Q(c_{1})\o c_3+ Q(c_{2})\o Q(c_{3})\o c_1\\
 &&~~~~~~=0,
 \end{eqnarray*}
 finishing the proof.                        \hfill $\square$

 \smallskip

 {\bf Theorem 4.3} Let $(C, \D, R)$ be a Rota-Baxter coalgebra of weight $0$. Define the operation $\widetilde{\D}$ on $C$ by
 $$
 \widetilde{\D}(c)=Q(c_1)\o c_2-Q(c_2)\o c_1.
 $$
 Then $\widetilde{C}=(C,\widetilde{\D})$ is a pre-Lie coalgebra.

 {\bf Proof.} Same to the proof of Theorem 4.2.                     \hfill $\square$

 \section{Rota-Baxter bialgebras}   \hskip\parindent
 In this section, we will combine the Rota-Baxter algebra and Rota-Baxter coalgebra to Rota-Baxter bialgebra and also provide some examples.

 {\bf Definition 5.1} Let $\lambda$, $\gamma$ be elements in $K$ and $H$ a bialgebra (maybe without unit and counit). A triple $(H, P, Q)$ is called a Rota-Baxter bialgebra of weight $(\lambda, \gamma)$ if $(H, P)$ is a Rota-Baxter algebra of weight $\lambda$ and $(H, Q)$ is a Rota-Baxter coalgebra of weight $\gamma$, that is, $H$ satisfies the following conditions: for all $x, y\in H$,
 $$
 P(x)P(y)=P(xP(y))+P(P(x)y)+\lambda P(x y),
 $$
 $$
 (Q\o Q)\D=(id\o Q)\D Q+(Q\o id)\D Q+\gamma\D Q.
 $$

 \smallskip
 {\bf Example 5.2} Let $H$ be a Hopf algebra and $(A, \mu, \D)$ a bialgebra. Suppose that $(A, \mu)$ is a right $H$-Hopf module algebra and $(A, \D)$ is a right $H$-Hopf module coalgebra. Then $(A, P_{R}, P_{R})$ is a Rota-Baxter bialgebra of weight $(-1,-1)$.

{\bf Proof.} By using Theorem 3.2 and \cite[Theorem 2.4]{Jian1}, we can get the conclusion directly.                       \hfill $\square$

 \smallskip
 {\bf Remark} Similarly, the above result holds for the left $H$-Hopf module algebra and the left $H$-Hopf module coalgebra.

 \smallskip

 {\bf Example 5.3}  Let $C$ be a bialgebra and $H$ a Hopf algebra with the antipode $S$. Suppose that there are two bialgebra maps: $i: H\lr C$ and $\pi: C\lr H$ such that $\pi\circ i=id_{H}$, i.e., $C$ is a bialgebra with a projection (see \cite{Ra}). Set $\Pi=id_{C}\star(i\circ S\circ \pi)$, where $\star$ is the convolution product on End($C$). The right $H$-Hopf module structure are given by the following:
 $$
 c\cdot h=ci(h),
 $$
 $$
 \rho_{R}(c)=c_1\o \pi(c_2),
 $$
 for all $c\in C$ and $h\in H$. Then $(C, \Pi, \Pi)$ is a Rota-Baxter bialgebra of weight $(-1, -1)$.

 {\bf Proof.} It is straightforward by using Example 3.3 and \cite[Example 2.6]{Jian1}.                \hfill $\square$

  \smallskip

 {\bf Example 5.4} Let $\{x, y, z\}$ be a basis of a 3-dimensional vector space $H$ over $K$. We define multiplication and comultiplication over $H$:
 $$
 x^{2}=x,~~ y^{2}=y,~~ z^{2}=0
 $$
 $$
 xy=yx=y, ~~yz=z,
 $$
 $$
 xz=zx=z, ~~zy=0,
 $$
 $$
 \D(x)=x\o x, ~~\D(y)=y\o y, ~~\D(z)=z\o z.
 $$
 It is not hard to see that $H$ is a bialgebra with the above constructions. Let $P_1, P_2, Q$ be operators defined with respect to the basis $\{x, y, z\}$ by
 $$
 P_{1}(x)=a z,~~P_{2}(x)=-c x-c y, ~~Q(x)=d z,
 $$
 $$
 P_{1}(y)=b z, ~~~~P_{2}(y)=-c y,  ~~~~Q(y)=d z,
 $$
 $$
 P_{1}(z)=0, ~~~~P_{2}(z)=-c z,  ~~~~Q(z)=0,
 $$
 where $a, b, c, d\in K$. Then $(H, P_{1}, Q)$ is a Rota-Baxter bialgebra of weight $(0, d)$, and $(H, P_{2}, Q)$ is a Rota-Baxter bialgebra of weight $(c, d)$.

 \smallskip

 {\bf Acknowledgments}
 The authors are deeply indebted to the referee for his/her very useful suggestions and some improvements to the original manuscript. This work was partially supported by the Foundation for Key Teacher by the Henan Normal University and the IRTSTHN (No. 14IRTSTHN023).

 \end{document}